\documentclass[11pt]{article}

\usepackage[letterpaper]{geometry}
\usepackage{microtype}
\usepackage{times}
\usepackage{graphicx}
\usepackage[colorlinks=true,allcolors=blue]{hyperref}

\title{Mathematicians in the Age of AI\footnote{I am grateful to Johan Commelin, Sidharth Hariharan, Bryna Kra, Emily Riehl, and Akshay Venkatesh for comments, corrections, and suggestions. A minor revision was made on April 6, 2026.}}

\author{Jeremy Avigad}

\date{March 3, 2026}

\begin{document}

\maketitle

\begin{abstract}
 Recent developments show that AI can prove research-level theorems in mathematics, both formally and informally. This essay urges mathematicians to stay up-to-date with the technology, to consider the ways it will disrupt mathematical practice, and to respond appropriately to the challenges and opportunities we now face.
\end{abstract}

\section{Background}

This past year was not an easy one for the U.S.~National Science Foundation, but last fall it managed to launch the Institute for Computer-Aided Reasoning in Mathematics (ICARM), even while facing serious budgetary constraints. The institute's mission is to support the adoption of new reasoning technologies in mathematics, many of them involving neural and symbolic AI. As director, I have given several variants of a talk titled ``Mathematics in the Age of AI,'' in which I survey the new technologies, including formalization and proof assistants, symbolic AI and automated reasoning, and machine learning and neural networks. These technologies are distinct, but they overlap and interact in interesting ways.

To date, early adopters of formal methods have contributed to a formal library for mathematics, Mathlib, and several collaborative projects have used it to verify contemporary mathematical results, even before or in the absence of the usual peer review process. Automated reasoning tools such as SAT solvers have solved open problems in combinatorics, algebra, number theory, and discrete geometry. Machine learning techniques have been used to detect patterns in and relationships between computational data, leading to the discovery of new results verified by conventional proofs. They have also been used to find combinatorial objects of interest, including counterexamples to conjectures. Neural networks are also used to compute solutions to PDEs, identify phenomena of interest, and search for parameter settings that give rise to such phenomena. Although the technologies are still niche and have seen only a handful of notable successes, the advances make it clear that they are bound to have a significant impact on the subject.

Machine learning and symbolic methods come together in the use of neural networks to prove mathematical theorems. Sometimes this involves using language models and agentic systems to write informal mathematical arguments, and sometimes it involves using them to construct formal proofs whose correctness is verified by a proof-checker like Lean. Many systems interleave reasoning on both levels, leveraging insights from conventional mathematical literature and the strong signal for correctness provided by a formal checker. Recent advances in this area have been especially dramatic.

The tone of my surveys has generally been optimistic: the new technologies can help us do mathematics better and do things we couldn't do before. They can have positive effects not just on the verification and discovery of mathematical results, but also on how we teach and learn mathematics, curate mathematical knowledge, and collaborate and communicate mathematical ideas. The talks are, however, tempered by some of the concerns about AI. I have treated these in greater depth in an essay,\footnote{\url{https://arxiv.org/abs/2502.14874}} ``Is Mathematics Obsolete?'' and, more recently, in a talk titled ``How to Worry About Mathematics in the Age of AI.'' I find it helpful to distinguish between two types of concerns, namely, concerns about AI in general and concerns about how AI will impact mathematics in particular. The general concerns have to do not only with the social, economic, and environmental effects of AI, but also the effects that AI will have on our ability to think and reason for ourselves, and the wisdom of putting our lives and livelihoods in the hands of systems whose behaviors we can't control and don't understand. I have argued in ``Obsolete'' that mathematics is at least part of the solution, since it allows us to ask precise questions and request explanations and justification in the form of artifacts that we can check, query, and audit independently. It gives us ways to interact with AI that keep us at the center of our deliberative processes and preserve our agency over our decisions.

This essay is about concerns of the second type, specifically, some of the ways that AI will impact the mathematical profession. It is prompted by recent developments that should give us pause. It is becoming increasingly clear that the technology will be disruptive, and we are not taking that fact seriously enough. As the director of an institute meant to support the adoption of the new technologies, I may be complicit in the disruption, but the point of this essay is that AI is here whether we like it or not, and if we don't rise to the occasion, we are all complicit. We owe it to our students, and to mathematics, to get this right.

\section{Recent Developments}

In early 2024, Sidharth Hariharan, then an exchange student from Imperial College visiting EPFL and now a first-year PhD student at Carnegie Mellon, met with Maryna Viazovska and began a project to formalize Viazovska's proof of the optimality of the $E_8$ lattice for 8-dimensional sphere packing. The approach they took is by now standard in the formalization community: they wrote a detailed informal blueprint to guide the formalization. Chris Birkbeck and Seewoo Lee joined the project over the summer, and in the fall, it became the basis of Sid's master's thesis, supervised by Bhavik Mehta. Viazovska gave a progress report in a talk at the Big Proof conference at the Isaac Newton Institute in Cambridge, UK, in July 2025, at which point the team opened the project to the public.

In the fall of 2025, things were going swimmingly. The project leadership used social media and weekly Zoom meetings to coordinate a team of volunteers. Several AI startup companies approached them seeking problems to test their products. The team embraced the use of AI, accepting AI-generated ``pull requests,'' revising them to meet their contribution standards, and incorporating them into the codebase. The project leaders enjoyed experimenting with AI themselves. Then, a few weeks ago, Sid showed up at my office door, looking pale; he had just learned that a company known as Math Inc.~decided to show off the abilities of its newest proving agent, Gauss, by completing a formal proof of the final result. The company had contributed to the project collaboratively in November, but then suddenly went dark, choosing to work in secret. This was \emph{not} the type of AI collaboration the group had in mind; they had not expected the company to throw computational resources at the project to obtain a splashy result, an act which Matt Ballard has aptly described as a ``drive-by proving.'' The company has since gone on to use Gauss to verify the corresponding 24-dimensional result, which builds on the 8-dimensional result, without additional input from the team.\footnote{See the announcement at \url{https://www.math.inc/sphere-packing}.}

Several prominent members of the Lean community, including the team of maintainers, rallied to help the project leaders respond to the news. One concern was the way the company would frame the result in its press release: Gauss's success was built on almost two years of creative work, scaffolding, and planning by the human participants, and it would have been unfair to them, mostly early-career researchers, to advertise this as solely a success for AI. A bigger concern was that the company would proclaim the project ``done.'' The formalization, on its own, is close to worthless, since the correctness of Viazovska's result was never in doubt. The participants embraced the project, rather, as a way of revisiting those results and better understanding them, and of building libraries and infrastructure to support future work. AI-generated proofs tend to be verbose and narrowly focused on the task at hand, and project leaders worried that a selfish credit grab would eliminate any incentive to revise the code and get it where they wanted it to be.

Math Inc.~has, however, cooperated with the project leaders to ensure that the announcements provide appropriate credit to the team and has promised to collaborate with the team openly to revise and improve the code to meet project standards. The community is now cautiously optimistic that this model for human-AI collaboration will work out well. A negative framing of the outcome is that the project organizers are making the best of a bad situation, while the positive framing is that we are all learning how to put AI to good mathematical use. It is still too early to say how this will play out; if AI can alleviate some of the tedium in formalization and help us enjoy and understand the mathematics, we will be better off, but if it gets in the way of that enjoyment and understanding, we will have lost something important.

While this episode in formal theorem proving was playing out, the mathematics community witnessed a complementary experiment in informal theorem proving. On February 5, a group of eleven mathematicians posted ten challenge problems for AI on arXiv.\footnote{\url{https://arxiv.org/abs/2602.05192}} Their goal was to assess the ability of AI to assist with real mathematical research. To that end, the group chose questions they had encountered in their own research and had answered with proofs that were about five pages long but not yet published. They tested two commercially available systems, Gemini 3.0 Deep Think and ChatGPT 5.2 Pro, and the results were not impressive. However, they also invited AI developers to try their own systems on the problems and make the results public before February 13, so the community could assess them as well.

A few AI developers took the bait, and the authors of the paper asked ICARM to host a discussion of the results in its Zulip social media platform,\footnote{\url{http://icarm.zulipchat.com/}, in the ``first proof'' channel} and then invited mathematicians to evaluate the AI contributions. A model produced by OpenAI offered a solution to one problem that experts judged to be ``completely correct and also quite beautiful.'' The best performance was by a proving agent, Aletheia, developed by Google DeepMind, which produced correct results for six out of the ten problems.\footnote{See the announcement, \url{https://arxiv.org/pdf/2602.21201}.} I suspect that I was not alone in being surprised by this level of success.

These developments show that AI is getting good at both informal and formal mathematics. When ChatGPT was first released in November 2022, we all laughed at how bad it was at answering mathematical questions. Last summer, four companies claimed gold-medal performance on the International Mathematical Olympiad competition problems, two generating informal solutions and two generating formal ones. The best systems are now well on their way to saturating Putnam benchmarks. Now that AI has begun to handle proofs at the level of mathematical research, we are running out of places to hide. AI has gone from close to zero proving power to the current state of affairs in less than four years, and investment in AI research and development has exploded in that time. We have to face up to the fact that AI will soon be able to prove theorems better than we can.

\section{Concerns}

We are now in an uncomfortable position of seeing AI do things we once thought we were uniquely able to do. We have been proud of our ability to solve really hard problems, construct long arguments with precise reasoning, and juggle lots of complex parts and put them together just right, but now we know that AI can do that too. We have also been proud of our deep insights and intuitions, our ability to detect subtle patterns and analogies, and our capacity to formulate far-reaching generalizations. Contemporary uses of AI are moving in that direction as well.

We need to think about what this will do to our profession. Why do we academic mathematicians have jobs? In part, it is because there are students majoring in mathematics, that is, young people who enjoy doing mathematics and whose parents are willing to support them. We should wonder what will happen to that enjoyment if doing mathematics becomes a matter of querying AI. Colleagues often point out to me that we still have professional chess players, even though chess engines are much better than any human player, but I find that small comfort: we don't require schoolchildren to play chess from kindergarten through high school, and it's hard to make a living as a chess tutor. Comparing the enjoyment of mathematics to the enjoyment of music and art may be a better analogy, but those disciplines, too, struggle for support. If the best argument we can muster for teaching mathematics is that it provides aesthetic or recreational enjoyment, we are in for a rough time.

A more substantial reason academic mathematicians have jobs is that we provide service courses for engineering, business, computer science, data science, and other fields that are well supported by industry and government. It's easy to imagine that source of support diminishing as well. Students can already use AI to do the kind of homework we have been assigning for decades, and using in-class exams to discourage that artificially deprives them of relying on tools that they will have ready access to in their careers. If we are not training engineers to use AI to solve problems, we are not training them to be engineers, and engineering departments won't need our courses. We might not expect demand to vanish overnight, but the academic environment is changing rapidly, and when the bottom falls out, it will be hard to recover.

\section{Recommendations}

Rather than let the dire alarms of the previous section bring us down, we should have faith in the power of mathematics.  As I have argued in ``Obsolete,'' the use of precise language and abstraction to help us understand the world around us, reason efficiently, and communicate ideas is as important as ever in the age of AI. Advances in technology provide us with new ways of achieving our mathematical goals, but AI is nothing more than technology, and it's up to us to decide how to use it.

At the end of ``Obsolete,'' I offered a dystopian view in which we turn our back on mathematics and cede agency over our deliberative processes to AI, and an optimistic view, in which mathematics remains central to those processes. The dystopian view translates roughly to the situation described in the previous section. How does the more positive view play out for mathematicians? The answer is simple: it has us using AI to do the kind of mathematics we know and love, but better. If AI can help us realize the Langlands program, prove Goldbach's conjecture, and resolve the $P \ne NP$ problem, will that be all that bad? We should look forward to the benefits of AI in both pure and applied mathematics. If scientists, engineers, and analysts can design more efficient systems and obtain more reliable results because mathematics provides them with the qualitative and quantitative tools to use AI effectively, mathematics will be doing more than enough to earn its keep.

We need to remember our strengths: mathematicians are problem solvers and theory builders \emph{extraordinaire}. Rather than fight the use of AI in mathematics, we should \emph{own} it. It is not enough to keep up with current events and design benchmarks for AI researchers; we need to play an active role in deploying the technology and molding it to our purposes. We also need to learn how to raise our students with the wisdom to use the new technologies appropriately, and we need to be careful that we still manage to impart core mathematical intuitions and understanding. Figuring out how to use AI effectively to achieve our mathematical goals won't be easy, but mathematicians have always embraced challenges---indeed, the harder, the better. If we face AI head-on and stay true to our values, mathematics will thrive. We just need to show up and get to work.

\end{document}